\documentstyle[amsfonts, 11pt]{article}

\newtheorem{backgroundtheorem}{Background Theorem}
\newtheorem{definition}{Definition}
\newtheorem{theorem}{Theorem}

\newtheorem{proposition}{Proposition}

\title{\sc Trace Functionals of the Kontsevich Quantization}

\author{Alexander Golubev\\ 
Max-Planck-Institute for Mathematics\\
Vivatsgasse 7, 53111 Bonn, Germany\\
\textit{e-mail address:}
\texttt{golubev@mpim-bonn.mpg.de}}
\input psfig.sty

\begin{document}

\maketitle

{\bf Abstract:}  {\small We generalize the notion of trace to the Kontsevich quantization algebra and show that for all Poisson manifolds representable by quotients of a symplectic manifold by a Hamiltonian action of a nilpotent Lie group, the trace is given by integration with respect to a unimodular volume form.}

\section{Introduction}

The gist of deformation quantization of a Poisson manifold lies in modifying the
usual (commutative) product of functions to obtain a new operation (star product) in such a way that the first approximation is given by the Poisson bracket. This
approach was first proposed in \cite{BFFLS}, and is at the moment one of the most
developed methods to relate "classical" data represented by a Poisson manifold to a 
"quantum" associative algebra. With the advent of the Kontsevich star product in 1997
(see \cite{Ko}), when the problem of existence of these algebras for arbitrary
Poisson structures has been solved in full generality, it became abundantly clear that
the utility of deformation quantization hinges on whether the resulting algebras are amenable to the standard operator theory tools, notwithstanding the fact that
no actual operators are present. If this universal version of quantization is to be useful, some "spectral"
properties must transpire. Sure enough, on the infinitesimal level there is the notion of Poisson trace, whose existence is intimately linked with the modular class (or, in mundane terms, with the existence of a unimodular measure - one invariant with respect to
all Hamiltonian flows) as demonstrated by Weinstein in
\cite{W2}. Another bit of information supporting the idea of the Kontsevich quantization algebras having spectral properties is the fact that on ${\mathbb{R}}^d$ any constant coefficient Poisson structure induces a Kontsevich quantized algebra with a trace given by integration.\\ 
\indent
On the other hand, considering all possible star products
on symplectic manifolds, Connes et al. \cite{CFS} were able to identify obstructions
to the existence of a trace on a fixed star-product algebra. These turned out to be
certain cocycles in the cyclic cohomology of the Hochshild complex. Moreover, the
authors constructed a simple example where integration does not yield a trace
functional. Combined with the fact that the Liouville measure on a symplectic
manifold is unimodular \textit{par excellence}, the results of Connes et al. permit one to conclude that unimodularity
alone is not sufficient even in the nicest possible (symplectic) setting. Although 
later Fedosov \cite{F} showed that on an arbitrary symplectic manifold there is
a star-product quantization equipped with a trace given by integration, the following
question remains: Is there a trace functional on the Kontsevich quantization algebra
of an arbitrary Poisson manifold? The present paper aims to single out a class of 
Poisson manifolds that possess a trace functional on the Kontsevich quantization 
algebra. In the process we first briefly recall the essentials in Section 2, then
study the Kontsevich quantization of the dual of nilpotent Lie algebras endowed with the standard Lie-Poisson structure (Section 3), go on to describe the properties of Morita equivalent Poisson manifolds in Section 4, and finally get everything together in Section 5.\\
\noindent
{\bf Aknowledgments.} We would like to thank Joseph Donin for many valuable suggestions, and Alan Weinstein for useful comments and a careful reading of the manuscript.

\section{Star Products}
Given a smooth (${\cal C}^{\infty}$) Poisson manifold $(P, \pi)$, the set of smooth
functions ${\cal C}^{\infty}(P)$ can be viewed as a commutative algebra. The star
product on ${\cal C}^{\infty}(P)$ (c.f. \cite{BFFLS}) is an associative ${\mathbb R}[[\hbar]]$-linear product on ${\cal C}^{\infty}(P)[[\hbar]]={\mathbb{Q}}(P)$ expressed by the following formula for $f, g \in {\cal C}^{\infty}(P) \subset {\mathbb{Q}}(P)$:\\
$$(f,g) \mapsto f\;\star\;g\;=\;fg\;+\;\hbar\pi(f,g)\;+\;{\hbar}^2B^2(f,g)\;+...
\in {\mathbb{Q}}(P),$$
where $\hbar$ is a formal variable, and $B^i$ are bidifferential operators (bilinear
maps ${\mathbb{Q}}(P) \times {\mathbb{Q}}(P) \longrightarrow {\mathbb{Q}}(P)$ which are
differential operators of globally bounded order with respect to each argument). The 
product of arbitrary elements of ${\mathbb{Q}}(P)$ is defined by the condition of 
linearity over  ${\mathbb R}$ and $\hbar$-adic continuity:\\
$$\large(\sum_{n \geq 0} {\hbar}^n f_n \large)\star \large(\sum_{n \geq 0} {\hbar}^n g_n \large)\;:=\;\sum_{k, l \geq 0} {\hbar}^{k+l}f_k g_l\;+\;\sum_{k, l \geq 0,\;m \geq 1} {\hbar}^{k+l+m}B^m(f_k,g_l).$$
\indent
Now we proceed to give a brief account of the universal deformation quantization of ${\mathbb R}^d$ or a domain thereof due to Kontsevich \cite{Ko}. From this point on,
the symbol $\star$ refers to the canonical Kontsevich star product.
In order to describe the terms proportional to ${\hbar}^n$ for any integer $n \geq 1$, Kontsevich introduced a special class $G_n$ of oriented labeled graphs called {\textit{admissible graphs}}.
\begin{definition}
An oriented graph $\Gamma$ is admissible ($\Gamma \in G_n$, $(n \geq 1)$) if:
\begin{enumerate}
\item $\Gamma$ has $n+2$ vertices labeled $\{1,2,...n,L,R\}$ where L and R stand for
Left and Right respectively, and $\Gamma$ has $2n$ oriented edges labeled $\{i_1,j_1,
i_2,j_2,...,i_n,j_n\}$;
\item The pair of edges $\{i_m,j_m\}, 1\leq m \leq n$ starts at the vertex m;
\item $\Gamma$ has no loops (edges starting at some vertex and ending at that same vertex)
and no parallel multiple edges (edges sharing the same starting and ending vertices).
\end{enumerate}
\end{definition}
\noindent
The class $G_n$ is finite. For $n \geq 1$ the first edge $i_k$ starting at the vertex $k$
has $n+1$ possible ending vertices since there are no loops, while the second edge $j_k$ has only $n$ possible "landing sites" since there is no parallel multiple edges.
Thus there are $n(n+1)$ ways to draw the pair of edges starting at some vertex and 
therefore $G_n$ has ${n(n+1)}^n$ elements. For $n=0$, $G_0$ has only one element:
the graph having $\{L,\;R\}$ as set of vertices and no edges. Of all admissible graphs, 
there is an important particular subclass, which we call {\textit{the class of unions of subgraphs}}. This notion will be used in the sequel. Here is how we define the elements
of this subclass:
\begin{definition}
A graph $\Gamma \in G_r$ is the union of two subgraphs ${\Gamma}_1$ and ${\Gamma}_2$ with
$r_1+r_2=r$, if the subset $(1,...,r)$ of the set of vertices of  $\Gamma$ can be split
into two parts $(a_1,...,a_{r_1})$ and $(b_1,...,b_{r_2})$ such that there is no edge 
between these two subsets of vertices.
\end{definition}
\indent
A bidifferential operator $(f,g)\mapsto B_{\Gamma}(f,g),\;f,g \in {\cal C}^{\infty}
({\mathbb R}^d)$ is associated to each graph  $\Gamma \in G_n,\;n\geq 1$. To each vertex
$k$, $1 \geq k \geq n$ one associates the components ${\pi}_{i_k j_k}$ of the Poisson
structure, $f$ is associated to the vertex $L$ and $g$ to the vertex $R$. Each edge 
such as $i_k$ acts by partial differentiation with respect to $x_{i_k}$ on its ending
vertex. For $n=0$ we simply have the usual product of $f$ and $g$. In what follows we use the words "graph" and "bidifferential operator encoded by the graph" interchangeably.\\
\indent
Now we go on to describe the weights $w_{\Gamma}$. Let ${\cal H} \;=\;\{ z \in {\mathbb C}|
{\textnormal{Im}} (z) > 0 \}$ be the upper half-plane. ${\cal H}_n$ will denote the
configuration space $\{z_1,...,z_n \in {\cal H}|z_i \ne z_j\;{\textnormal{for}}\;i \ne j\}$.
${\cal H}_n$ is an open submanifold of ${\mathbb C}^n$. Let ${\phi}:{\cal H}_2 \rightarrow
{\mathbb{R}}/2{\pi}{\mathbb{Z}}$ be the function:\\
$${\phi}(z_1,z_2)\;=\;{1\over{2 \sqrt{-1}}}{\textnormal{Log}}{(z_2 -z_1)({\bar{z}}_2-z_1)
\over{(z_2 - {\bar{z}}_1)({\bar{z}}_2-{\bar{z}}_1)}}.$$
\noindent
${\phi}(z_1,z_2)$ is extended by continuity for $z_1,\;z_2 \in {\mathbb{R}},\;z_1 \ne z_2$.\\
\indent
For a graph $\Gamma \in G_n$, the vertex $k$, $1 \leq k \leq n$ is associated with
the variable $z_k \in {\cal H}$, the vertex $L$ with $0 \in {\mathbb{R}}$, and the vertex
$R$ with $1 \in {\mathbb{R}}$.\\
\indent
The weight $w_{\Gamma}$ is defined by integrating an $2n$-form over ${\cal H}_n$:\\
$$w_{\Gamma}\;=\;{1 \over{n! {(2 \pi)}^{2n}}} \int_{{\cal H}_n} \bigwedge_{1 \leq k \leq n}
(d{\phi}(z_k,\;I_k) \wedge d{\phi}(z_k,\;J_k)),$$
where $I_k$ ($J_k$) denotes the variable or real number associated with the ending vertex of 
the edge $i_k$ ($j_k$). It is showed in \cite{Ko} that this integral is absolutely convergent.
As one can see, the weights do not depend on the Poisson structure or
the dimension of the underlying manifold. Combining these constructions Kontsevich proved
\begin{backgroundtheorem}[Kontsevich] \hspace{.01in}
Let $\pi$ be a Poisson structure in a domain of ${\mathbb{R}}^d$. The formula\\
$$f\; \star \;g := \sum_{n=0}^{\infty} {\hbar}^n \sum_{\Gamma \in G_n}
w_{\Gamma} B_{\Gamma, \pi} (f, g)$$
defines an associative product up to gauge equivalence.
\end{backgroundtheorem}
\noindent
Furthermore, all the machinery above as well as 
Background Theorem 1 generalize to arbitrary Poisson manifolds using localization. There are no cohomological obstructions.\\
\indent
Finally, we deal with trace functionals on star-product algebras.
Consider the linear map  ${\mathsf Tr}$ on the compactly supported Kontsevich quantization algebra ${\mathbb{Q}}_c(P)$
with values in formal Laurent power series ${\hbar}^{-[d/2]}{\mathbb R}[[\hbar]]$:
$${\mathsf Tr}_{P,{\mu}}(F)\;=\;{\hbar}^{-[d/2]}\int_{P}F{\mu},$$
where [ ] denotes the integral part of a number and $d$ is the dimension of $P$. When
$P$ is symplectic this functional (up to a contant) coincides with the usual one as 
defined in \cite{CFS}, provided $\mu$ is the Liouville volume form. In general, we
may have many meaningful choices for the volume form, which is emphasized by the subscript.
Following Fedosov \cite{F} we declare that the functional ${\mathsf Tr}$ satisfies
the trace property with respect to the Kontsevich star product if
$${\mathsf Tr}(F \star G)\;=\;{\mathsf Tr}(G \star F).$$
\noindent
Consequently, the problem of defining a trace functional amounts to proving that 
a linear functional ${\mathsf Tr}$ satisfies the trace property. Formally, we have
\begin{definition}
The linear functional ${\mathsf Tr}$ on ${\mathbb{Q}}_c(P)$ is called a trace functional
if it satisfies the trace property. 
\end{definition}

\section{Quantization of ${\mathfrak{g}}^*$}
Given a finite-dimensional Lie group $G$, we denote its Lie algebra by $\mathfrak g$,
the dual of $\mathfrak g$ by $\mathfrak g^*$. The structure constants of $\mathfrak g$
are determined by the set of relations with respect to a basis of  $\mathfrak g$:\\
$$[X_i, X_j]\;=\;X_k C^k_{ij},\;\; X_i \in  {\mathfrak g}\;\;\forall i \in \{1, ...., d\}.$$
\noindent
Now we let $(x_1, ...., x_d)$ denote linear coordinates on  $\mathfrak g^*$. Then 
the natural Lie-Poisson structure\\
$$\pi \; = \sum_{i, j} \pi_{ij} \partial_i \wedge \partial_j,$$
is expressed in terms of coordinates via
$$\pi_{ij} \; =\; x_k C_{ij}^k .$$
\indent
We begin by looking at the role of unimodularity in vanishing of the Poisson trace. 
\begin{proposition} \hspace{.01in}
For any $f,\;g \in {\cal{C}}_c^{\infty}(\mathfrak g^*)$, \\
$$ \int_{{\mathbb R}^d} \pi(f,g) dx_1 \wedge ...\wedge dx_d\;=\;0$$
if and only if  $$\sum_{i}  C_{ij}^i\;=\;0\;\;\forall j \in \{1, ...., d\}.$$
\end{proposition}
\noindent
{\textit {Proof.}}  We proceed by making a simple observation: $C^k_{ij} {\partial}_i
{\partial}_jf$ adds up to zero by virtue of $C^k_{ij}$ being antisymmetric. Now straightforward integration by parts of a single component $\pi_{ij}{{\partial}_i} f{{\partial}_j}g$ with subsequent summation does the trick. Q.E.D.\\
\noindent
{\textit{Remark.}}
Proposition 1 is basically a rephrasing of the statement made in \cite{W2}, according to which "on the dual of a Lie algebra $\mathfrak g$, with its Lie-Poisson structure, the modular vector field with respect to any translation-invariant density is the constant vector field with value ${\mathsf{Tr}}\;ad$, the trace of the adjoint representation".\\
\indent
Now we are in a position to state the main theorem of this section.
\begin{theorem} \hspace{.01in}
For all $F, G\;\in {\mathbb Q}_c({\mathfrak g}^*)$ the equality\\
$$ \int_{{\mathbb R}^d} F \star G dx_1 \wedge ...\wedge dx_d\;=\;
\int_{{\mathbb R}^d} G \star F dx_1 \wedge ...\wedge dx_d$$
holds for $\mathfrak g$ nilpotent.
\end{theorem}
\noindent
{\textit{Proof.}}  First of all we note that this Theorem will follow if we manage
to show the equality above for $f, g \in {\cal C}_c^{\infty}({\mathfrak g}^*)$.
Next we invoke the formula (due to Cattaneo and Felder \cite{CaFe})
relating the weight of $B^n(f,g)$ to that of  $B^n(g,f)$. Namely, $w_{B^n(f,g)}=
(-1)^nw_{B^n(g,f)}$. This signifies that for all even $n$ the corresponding integrals
do not affect the (presumed) trace property. From now on we work only with the graphs involving odd $n$. To manipulate the operators effectively, we introduce the 
following notation: $B^n_k(f,g)$ stands for a bidifferential operator that  has order
$n$ in the star product formula, and, in addition, has exactly $k$ edges acting 
summarily on $f$
AND $g$. Nothing else is specified, but below we take one particular operator at a
time so that this apparent ambiguity does not matter. Now, by the linearity of $\pi$,
it follows that $n\leq k \leq 2n$. An important fact used later is that any graph
denoted by $B_n^n(f,g)$ necessarily includes at least one loop (an expression of the form
$\partial_{i_m}{\pi}_{{i_1}{j_1}}\partial_{i_1}{\pi}_{{i_2}{j_2}}...\partial_{i_{m-1}}
{\pi}_{{i_m}{j_m}}$ after possible reshuffling). As Kathotia \cite{Ka} showed, such graphs vanish for nilpotent Lie algebras due to a theorem of Varadarajan \cite{V}:
\begin{backgroundtheorem}[Varadarajan, Theorem 3.5.4.] \hspace{.01in}
Let ${\mathfrak g}$ be a d-dimensional nilpotent Lie algebra over a field of 
characteristic zero. Then there is a basis $\{X_1, ..., X_d\}$ of ${\mathfrak g}$
such that the structure constants $C^k_{ij}$ defined by $[X_i,X_j]=X_kC^k_{ij}$ satisfy
$$C^r_{ij}\;=\;0,\;\;\forall r\geq min(i, j).$$
\end{backgroundtheorem}
An easy corollary to the above theorem assures that any cyclic product of the form
$C_{i_1j_1}^{j_m}C_{i_2j_2}^{j_1}...C_{i_mj_m}^{j_{m-1}}\;(m>1)$ vanishes and this,
in turn, yields the vanishing of loop graphs.\\
\indent
Now since we have seen that in the nilpotent setting there are no loop graphs, the idea is  to rewrite the integrals of $B_k^n$ as combinations of those of $B_n^n$. This is
accomplished via integration by parts. Once again using the linearity of the Poisson
structure, we have a vertex such that one outgoing edge acts on $f$, whereas the other
acts on $g$ (otherwise $k=n$ and we are done). We use one of these edges to perform 
integration by parts. As a result we obtain a few integrals
of admissible graphs of the form $B_{k-1}^n$. To see that we first note that all nilpotent groups are unimodular. Indeed, if there is a $C_{ij}^i \ne 0$ with respect to any basis, we would have an $X_i$ such that $[X_i,X_j]=X_iC^i_{ij}$, and iterating this $l$ times we would violate the condition ${\mathfrak{g}}^l=0$. So the admissibility of the resulting 
graphs as well as the decreasing of the number of edges acting on $f$ and $g$ follows 
from Proposition 1 and its proof. Inductively we get rid of all antisymmetric ($n$ odd) graphs. Q.E.D.\\
\noindent
{\textit{Remark.}} This result is in some sense the best attainable. Consider a unimodular Lie algebra ${\mathfrak u}$ such that the Kontsevich quantization of  ${\mathfrak u}^*$
does contain loop graphs. Then by choosing $f, g$ supported in a small neighborhood
of the origin in ${\mathfrak u}^*$, we can assure that the contribution of the graphs
with $n$ odd, $k=n$, is not eliminated by some freaky cancellation. Thus the statement
of Theorem 1 is not true for arbitrary unimodular Lie algebras.\\

\section{Morita Equivalence}
The essential sources used here are \cite{W1}, \cite{X}, and \cite{GG}. 
Following \cite{W1} we recall that a {\textit{full dual pair}}  $P_1 \stackrel{{\rho}_1}{\leftarrow}W \stackrel{{\rho}_2}{\rightarrow}P_2$
consists of two Poisson manifolds $(P_1, {\pi}_1)$ and $(P_2, {\pi}_2)$, a 
symplectic manifold $W$, and two submersions ${\rho}_1: W \rightarrow P_1$ and
${\rho}_2: W \rightarrow P_2$ such that ${\rho}_1$ is Poisson, ${\rho}_2$ is 
anti-Poisson, and the fibers of ${\rho}_1$ and ${\rho}_2$ are symplectic orthogonal
to each other. A Poisson (or anti-Poisson) mapping is said to be  {\textit{complete}}
if the pullback of a complete Hamiltonian flow under this mapping is complete. A full
dual pair is called complete if both ${\rho}_1$ and ${\rho}_2$ are complete. The
Poisson manifolds $(P_1, {\pi}_1)$ and $(P_2, {\pi}_2)$ are {\textit{Morita equivalent}} if there exists a complete full dual pair  $P_1 \stackrel{{\rho}_1}{\leftarrow}W \stackrel{{\rho}_2}{\rightarrow}P_2$ such that
${\rho}_1$ and ${\rho}_2$ both have connected and simply connected fibers. The notion
of Morita equivalence of Poisson manifolds was introduced and studied by Xu \cite{X},
as a classical analogue of the Morita equivalence of $C^*$-algebras.

\begin{theorem} \hspace{.01in}
The property of being equipped with a trace functional on the Kontsevich
quantization algebra is an invariant of Morita equivalence.
\end{theorem}
\noindent
{\textit{Proof.}}
By the hypothesis we have $P_1$, $P_2$ Morita equivalent, and ${\mathbb{Q}}_c(P_1)$
possesses a trace functional. Whence\\
$${\mathsf Tr}_{P_1,{\mu}'}(F)\;=\;{\hbar}^{-[d/2]}\int_{P_1}F{\mu}',$$
where ${\mu}'$ is a unimodular volume form. Now we recall a result from \cite{GG}. Here
mod$(P_i)$ denotes the modular class on respective manifolds.
\begin{backgroundtheorem}[Ginzburg, Golubev] \hspace{.01in}
Let $P_1$ and $P_2$ be Morita equivalent and let, in addition, $P_1$ be locally
unimodular. Then $P_2$ is also locally unimodular and $\em{mod}(P_1)$ goes to
$\em{mod}(P_2)$ under the natural isomorphism of the first Poisson cohomology groups
$E: H_{\pi}^1(P_1) \stackrel{\cong}{\rightarrow}H_{\pi}^1(P_2)$, i. e.
$E({\em{mod}}(P_1))={\em{mod}}(P_2)$.
\end{backgroundtheorem}
\noindent
Using the above theorem we conclude $P_2$ is unimodular too. Furthermore, starting off with ${\mu}'$, utilizing the action of Brylinski's symplectic star operator (see 
\cite{Br} and also \cite{GG} for more details), we arrive at a global unimodular volume form ${\mu}''$ on $P_2$. At this point we invoke a theorem of Weinstein \cite{W1} concerning transversal Poisson structures in full dual pairs:
\begin{backgroundtheorem}[Weinstein] \hspace{.01in}
Let $P_1 \stackrel{{\rho}_1}{\leftarrow}W \stackrel{{\rho}_2}{\rightarrow}P_2$ be a full dual pair. For each $x\in W$, the transverse Poisson structures  on $P_1$ and $P_2$ at ${{\rho}_1}(x)$ and ${{\rho}_2}(x)$ are anti-isomorphic. Consequently, if ${\em dim}(P_1) \leq {\em dim}(P_2$), $P_2$ is locally anti-isomorphic to the product of $P_1$ with a symplectic manifold.
\end{backgroundtheorem}
\noindent
Now we introduce a fixed open cover $\{U_{\alpha}\}_{\alpha \in J}$ of $W$ such that
$U'_{\alpha}={\rho}_1(U_{\alpha})$, $U''_{\alpha}={\rho}_2(U_{\alpha})$
and the latter ones are open covers of $P_1$ and $P_2$ respectively. They are so
fine that there is a local anti-isomorphism of Background Theorem 4 in each $U_{\alpha}$ 
(by refining the cover of $W$ we can always achieve that).
From Background Theorem 4 we infer the relation between the local expressions for
${\mu}'$ and ${\mu}''$, which we denote by ${\mu}'_{\alpha}$ and ${\mu}''_{\alpha}$.
Namely, either ${\mu}'_{\alpha}={\mu}''_{\alpha}\wedge \kappa$ (this happens to be the 
case if ${\textnormal{dim}}(P_1) > {\textnormal{dim}}(P_2)$), or  ${\mu}''_{\alpha}={\mu}'_{\alpha}\wedge \kappa$ (${\textnormal{dim}}(P_1) < {\textnormal{dim}}(P_2)$), or ${\mu}'=a{\mu}''$ (${\textnormal{dim}}(P_1) = {\textnormal{dim}}(P_2)$), where $\kappa$ is a Liouville volume form and $a$ is a nonzero constant. By the standard manifold theory we may view $U_{\alpha}$'s as domains of ${\mathbb 
R}^d$. A crucial fact needed at this juncture is the following: the graphs in $U'_{\alpha}$, $U''_{\alpha}$ are unions of subgraphs (c. f. Definition 2) of symplectic and transversal components. Moreover, the symplectic components necessarily have equal number of edges acting on $f$ and $g$, so that up to a constant, the graphs can be thought of as those involving transversal components only. Now assembling the above volume forms and
the graphs involving transversal components in the integral, we see that the vanishing of one on $U'_{\alpha}$ entails the vanishing of its counterpart on $U''_{\alpha}$.
Finally, using a partition of unity subordinate to the cover, we globalize the trace
functional. Q.E.D.

\section{The Main Theorem via Symplectic Realizations}
We recall the appropriate definitions.
A  {\textit{symplectic realization}} of a Poisson manifold $P$ is a pair $(W, \rho)$,
where $W$ is a symplectic manifold and $\rho$ is a Poisson morphism from $W$ to $P$. A
symplectic realization $\rho: W \rightarrow P$ is called  {\textit{complete}} if $\rho$ is complete as a Poisson map and $\rho$ is said to be  {\textit{full}} if it is a submersion. There are many ways to construct symplectic realizations, but we are interested in just one particular way of doing it. Precisely, we want to realize Poisson manifolds as quotients of symplectic manifolds by group actions. In applications, as Weinstein \cite{W1} pointed out, the symplectic manifold $W$ may represent a collection of states such that the points belonging to the same orbit of $G$-action are considered to be physically indistinguishable. Thus the set of "true physical states" is none other than the Poisson manifold $W/G$, and the group $G$ is called a  {\textit{gauge group}}. This motivates the definition below.
\begin{definition}
A Poisson manifold $P$ is symplectically realizable  with a gauge group $G$ if
there exists a symplectic manifold $W$ and a Lie group $G$ freely acting on $W$ by
symplectomorphisms, such that $W/G \cong P$. 
\end{definition}
\noindent
An upshot of this definition is that a) if $P$ is symplectically realizable with a gauge
group $G$, $\rho: W \rightarrow P$ is a quotient map and necessarily a complete full
surjective symplectic realization; b) $P$ is Morita equivalent to an open subset of ${\mathfrak g}^*$. 
\begin{theorem} \hspace{.01in}
The Kontsevich quantization algebra of a symplectically realizable
Poisson manifold with a nilpotent gauge group possesses a trace functional.
\end{theorem}
{\textit{Proof.}} Applying Theorem 1 we obtain a trace functional on ${\mathfrak g}^*$
in view of the hypothesis, and Theorem 2 now ensures that
${\mathsf Tr}_{P,\mu}$ satisfies the trace property, where $\mu$ is the volume form
obtained from the standard translation-invariant form on ${\mathfrak g}^*$. Q.E.D.\\

\end{document}